# COINCIDENCE AND COMMON FIXED POINT RESULTS FOR CONTRACTION TYPE MAPS IN PARTIALLY ORDERED METRIC SPACES

HASSEN AYDI

ABSTRACT. We present coincidence and common fixed point results of self-mappings satisfying a contraction type in partially ordered metric spaces. As an application, we give an existence theorem for a common solution of integral equations.

## 1. INTRODUCTION

Existence of fixed or common fixed results in partially ordered sets has fascinated many researchers in the few recent years where some applications to matrix equation, ordinary differential equations and integral equations are presented, we cite for example [1, 3, 4, 9, 13, 14].
Like in [9], we denote $\mathcal{S}$ the class of the functions $\beta : [0, \infty[ \to [0, 1[$ which satisfies the condition

$$(1.1) \qquad \beta(t_n) \to 1 \quad \text{implies} \quad t_n \to 0.$$

A. Amini-Harandi and H. Emami [9] proved recently in the contest of partially ordered sets the following result

**Theorem 1.1.** *Let $(X, \preceq)$ be a partially ordered set and suppose that there exists a metric $d$ in $X$ such that $(X, d)$ is a complete metric space. Let $T : X \to X$ be a non-decreasing mapping such that*

$$(1.2) \qquad d(Tx, Ty) \leq \beta(d(x,y))d(x,y), \quad \text{for} \ \ x, y \in X \ \ \text{with} \ \ x \preceq y,$$

*where $\beta \in \mathcal{S}$. Assume that either $T$ is continuous or $x$ satisfies the following condition : If $\{x_n\}$ is a nondecreasing sequence in $X$ such that $x_n \to x$, then $x_n \preceq x \ \forall n \in \mathbb{N}$.*
*Besides, suppose that for each $x, y \in X$, there exists $z \in X$ which is comparable to $x$ and $y$. If there exists $x_0 \in X$ with $x_0 \preceq Tx_0$, then $T$ has a unique fixed point.*

This previous theorem is an extension of Geraghty's result [7] where the contest of metric complete spaces is considered. This paper is again an extension of theorem 1.1. For partially ordered metric spaces, we take three self-mappings satisfying the following contraction

$$d(fx, gy) \leq \ \ \beta(d(Hx, Hy))\, d(Hx, Hy).$$

We present coincidence and common fixed point results of the self-mappings $f$, $g$ and $H$. These results will be given in section 2. As an application of the above results, in section 3, we give in a particular case an existence theorem for a common solution of integral equations.

## 2. Main results

We start some definitions, which we need in the sequel

**Definition 2.1.** *Let $X$ be a non-empty set, $N$ is a natural number such that $N \geq 2$ and $T_1, T_2, \cdots, T_N : X \to X$ are given self-mappings on $X$. If $w = T_1 x = T_2 x = \cdots = T_N x$ for some $x \in X$, then $x$ is called a coincidence point of $T_1, T_2, \cdots, T_{N-1}$ and $T_N$, and $w$ is called a point of coincidence of $T_1, T_2, \cdots, T_{N-1}$ and $T_N$. If $w = x$, then $x$ is called a common fixed point of $T_1, T_2, \cdots, T_{N-1}$ and $T_N$.*

**Definition 2.2.** [8] *Let $(X, d)$ be a metric space and $f, g : X \to X$ are given self-mappings on $X$. The pair $\{f, g\}$ is said to be compatible if $\lim_{n \to +\infty} d(fgx_n, gfx_n) = 0$, whenever $\{x_n\}$ is a sequence in $X$ such that $\lim_{n \to +\infty} fx_n = \lim_{n \to +\infty} gx_n = t$ for some $t$ in $X$.*

Let $X$ be a non-empty set and $R : X \to X$ be a given mapping. For every $x \in X$, we denote by $R^{-1}(x)$ the subset of $X$ defined by:

$$R^{-1}(x) := \{u \in X \mid Ru = x\}.$$

**Definition 2.3.** *Let $(X, \preceq)$ be a partially ordered set and $T, S, R : X \to X$ are given mappings such that $TX \subseteq RX$ and $SX \subseteq RX$. We say that $S$ and $T$ are weakly increasing with respect to $R$ if and only if for all $x \in X$, we have:*

$$Tx \preceq Sy, \; \forall \, y \in R^{-1}(Tx)$$

*and*

$$Sx \preceq Ty, \; \forall \, y \in R^{-1}(Sx).$$

**Remark 2.4.** *If $R : X \to X$ is the identity mapping ($Rx = x$ for all $x \in X$), then $S$ and $T$ are weakly increasing with respect to $R$ implies that $S$ and $T$ are weakly increasing mappings. Note that the notion of weakly increasing mappings was introduced in [2] (also see [5, 6]).*

Our first main result is an extension of theorem 1.1, and it is the following

**Theorem 2.5.** *Let $(X, \preceq)$ be a partially ordered set. Suppose that there exists a metric $d$ on $X$ such that $(X, d)$ is complete. Let $f, g, H : X \to X$ be given mappings satisfying*
*(a) $fX \subseteq HX$, $gX \subseteq HX$,*
*(b) $f$, $g$ and $H$ are continuous,*
*(c) the pairs $\{f, H\}$ and $\{g, H\}$ are compatible,*
*(g) $f$ and $g$ are weakly increasing with respect to $H$.*
*Suppose that for every $(x, y) \in X \times X$ such that $Hx$ and $Hy$ are comparable, we have*

(2.1) $$d(fx, gy) \leq \beta(d(Hx, Hy)) \, d(Hx, Hy),$$

*where $\beta \in \mathcal{S}$. Then, $f$, $g$ and $H$ have a coincidence point $u \in X$, that is, $fu = gu = Hu$.*

*Proof.* Let $x_0$ be an arbitrary point in $X$. Since $fX \subseteq HX$, there exists $x_1 \in X$ such that $Hx_1 = fx_0$. Since $gX \subseteq HX$, there exists $x_2 \in X$ such that $Hx_2 = gx_1$. Continuing this process, we can construct sequences $\{x_n\}$ and $\{y_n\}$ in $X$ defined by

(2.2) $$Hx_{2n+1} = fx_{2n} = y_{2n}, \quad Hx_{2n+2} = gx_{2n+1} = y_{2n+1}, \quad \forall \, n \in \mathbb{N}$$

By construction, we have $x_1 \in H^{-1}(fx_0)$ and $x_2 \in H^{-1}(gx_1)$, then using the fact that $f$ and $g$ are weakly increasing with respect to $H$, we obtain

$$Hx_1 = fx_0 \preceq gx_1 = Hx_2 \preceq fx_2 = Hx_3.$$

We continue the process to get

(2.3) $$Hx_1 \preceq Hx_2 \preceq ... \preceq Hx_{2n+1} \preceq Hx_{2n+2} \preceq ...$$

We can then write

(2.4) $$y_0 \preceq y_1 \preceq ... \preceq y_{2n} \preceq y_{2n+1} \preceq ...$$

<u>First case.</u> If there exists $n \in \mathbb{N}^*$ such that $y_{2n-1} = y_{2n}$, then by construction of the sequence $\{y_m\}$ and the contraction (2.1) with $x = x_{2n}$, $y = x_{2n+1}$.

$$\begin{aligned} d(y_{2n}, y_{2n+1}) = d(fx_{2n}, gx_{2n+1}) &\leq \beta(d(Hx_{2n}, Hx_{2n+1}))d(Hx_{2n}, Hx_{2n+1}) \\ &= \beta(d(y_{2n-1}, y_{2n}))d(y_{2n-1}, y_{2n}) = 0, \end{aligned}$$

which implies $y_{2n} = y_{2n+1}$. This leads to $y_m = y_{2n-1}$ for any $m \geq 2n$. Hence for every $m \geq 2n$ we have $Hx_m = Hx_{2n}$. This implies that $\{Hx_n\}$ is a Cauchy sequence. The same conclusion holds if $y_{2n} = y_{2n+1}$.

<u>Second case.</u> Suppose that $y_n \neq y_{n+1}$ for any integer $n$.
First, we will show that

(2.5) $$\lim_{n \to +\infty} d(y_{n+1}, y_n) = 0.$$

Thanks to (2.3), $Hx_{2n}$ and $Hx_{2n+1}$ are comparable, then using (2.2) and taking $x = x_{2n+2}$ and $y = x_{2n+1}$ in (2.1), we get

(2.6)
$$\begin{aligned} d(y_{2n+2}, y_{2n+1}) = d(Hx_{2n+3}, Hx_{2n+2}) &= d(fx_{2n+2}, gx_{2n+1}) \\ &\leq \beta(d(Hx_{2n+2}, Hx_{2n+1}))d(Hx_{2n+2}, Hx_{2n+1}) \\ &= \beta(d(y_{2n+1}, y_{2n}))d(y_{2n+1}, y_{2n}). \end{aligned}$$

Using $0 \leq \beta < 1$, we deduce then

(2.7) $$d(y_{2n+2}, y_{2n+1}) \leq d(y_{2n+1}, y_{2n}).$$

Similarly to this, one can find for $x = x_{2n}$ and $y = x_{2n+1}$ in (2.1) that

(2.8) $$d(y_{2n+1}, y_{2n}) \leq d(y_{2n}, y_{2n-1}).$$

Thus combining (2.7) together with (2.8) leads that for any $n \in \mathbb{N}$

(2.9) $$d(y_{n+2}, y_{n+1}) \leq d(y_{n+1}, y_n).$$

It follows that the sequence $\{d(y_{n+1}, y_n)\}$ is monotonic decreasing. Hence, there exists $r \geq 0$ such that

(2.10) $$\lim_{n \to +\infty} d(y_{n+1}, y_n) \to r.$$

From (2.6), we have

$$\frac{d(y_{2n+2}, y_{2n+1})}{d(y_{2n+1}, y_{2n})} \leq \beta(d(y_{2n+1}, y_{2n})) < 1.$$

Letting $n \to +\infty$ in the above inequality, then thanks to (2.10), we obtain

$$\lim_{n \to +\infty} \beta(d(y_{2n+1}, y_{2n})) = 1,$$

and since $\beta \in \mathcal{S}$, this implies that $r = 0$. Hence, (2.5) holds.
We need now to check that $\{Hx_n\}$ is a Cauchy sequence. Following (2.2), it suffices to prove that $\{Hx_{2n}\}$ is a Cauchy sequence. To do this, we proceed by contradiction. Suppose that $\{Hx_{2n}\}$ is not a Cauchy sequence. Then for any $\varepsilon > 0$, for which there exist two sequences of positive integers integers $\{m(k)\}$ and $\{n(k)\}$ such that for all positive integers $k$,

(2.11)
$$n(k) > m(k) > k, \quad d(Hx_{2m(k)}, Hx_{2n(k)}) > \varepsilon, \quad d(Hx_{2m(k)}, Hx_{2n(k)-2}) \leq \varepsilon$$

Therefore, we use (2.11) and the triangular inequality to get

$$\begin{aligned} \varepsilon &< d(Hx_{2m(k)}, Hx_{2n(k)}) \\ &\leq d(Hx_{2m(k)}, Hx_{2n(k)-2}) + d(Hx_{2n(k)-2}, Hx_{2n(k)-1}) + d(Hx_{2n(k)-1}, Hx_{2n(k)}) \\ &\leq \varepsilon + d(Hx_{2n(k)-2}, Hx_{2n(k)-1}) + d(Hx_{2n(k)-1}, Hx_{2n(k)}). \end{aligned}$$

Letting $k \to +\infty$ in the above inequality and using (2.5), we find

(2.12) $$\lim_{k \to +\infty} d(Hx_{2m(k)}, Hx_{2n(k)}) = \varepsilon.$$

Again, using the triangular inequality we have

$$\mid d(Hx_{2n(k)}, Hx_{2m(k)-1}) - d(Hx_{2n(k)}, Hx_{2m(k)}) \mid \leq d(Hx_{2m(k)}, Hx_{2m(k)-1}).$$

Letting again $k \to +\infty$ in the above inequality and using (2.5)-(2.12) we find

(2.13) $$\lim_{k \to +\infty} d(Hx_{2n(k)}, Hx_{2m(k)-1}) = \varepsilon.$$

On the other hand, we have

$$\begin{aligned} d(Hx_{2n(k)}, Hx_{2m(k)}) &\leq d(Hx_{2n(k)}, Hx_{2n(k)+1}) + d(Hx_{2n(k)+1}, Hx_{2m(k)}) \\ &= d(Hx_{2n(k)}, Hx_{2n(k)+1}) + d(fx_{2n(k)}, gx_{2m(k)-1}). \end{aligned}$$

Thanks to (2.5)-(2.12), then letting $k \to +\infty$, we have from the above inequality

(2.14) $$\varepsilon \leq \lim_{k \to +\infty} d(fx_{2n(k)}, gx_{2m(k)-1}).$$

We take now $x = x_{2n(k)}$ and $y = x_{2m(k)-1}$ in (2.1). Hence,

$$\begin{aligned} d(fx_{2n(k)}, gx_{2m(k)-1}) &\leq \beta(d(Hx_{2n(k)}, Hx_{2m(k)-1}))d(Hx_{2n(k)}, Hx_{2m(k)-1}) \\ &< d(Hx_{2n(k)}, Hx_{2m(k)-1}) \end{aligned}$$

Letting again $k \longrightarrow +\infty$ in the above inequality and using (2.5)-(2.13), we obtain

(2.15) $$\lim_{k \to +\infty} d(fx_{2n(k)}, gx_{2m(k)-1}) \leq \varepsilon.$$

Combining (2.14) to (2.15) yields

$$\lim_{k \to +\infty} d(fx_{2n(k)}, gx_{2m(k)-1}) = \varepsilon,$$

Therefore, since we are in the case $Hx_{2n(k)} \neq Hx_{2m(k)-1}$, then writing

$$\frac{d(fx_{2n(k)}, gx_{2m(k)-1})}{d(Hx_{2n(k)}, Hx_{2m(k)-1})} \leq \beta(d(Hx_{2n(k)}, Hx_{2m(k)-1})) < 1,$$

and using the fact that $\varepsilon = \lim_{k \to +\infty} d(fx_{2n(k)}, gx_{2m(k)-1}) = \lim_{k \to +\infty} d(Hx_{2n(k)}, Hx_{2m(k)-1})$, we get

$$\lim_{k \to +\infty} \beta(d(Hx_{2n(k)}, Hx_{2m(k)-1})) = 1.$$

We know that $\beta \in \mathcal{S}$, hence

$$\lim_{k \to +\infty} d(Hx_{2n(k)}, Hx_{2m(k)-1}) = 0,$$

which is a contradiction with (2.13), that is $\lim_{k \longrightarrow +\infty} d(Hx_{2n(k)}, Hx_{2m(k)-1}) = \varepsilon > 0$. We deduce then $\{Hx_n\}$ is a Cauchy sequence.

Let us now prove the existence of a coincidence point. First, we know that $\{Hx_n\}$

is a Cauchy sequence, then since $(X, d)$ is a complete metric space, there exists $u \in X$ such that

(2.16) $$\lim_{n \to +\infty} Hx_n = u.$$

From (2.16) and the continuity of $H$, we get

(2.17) $$\lim_{n \to +\infty} H(Hx_n) = Hu.$$

The triangular inequality yields

(2.18) $d(Hu, fu) \leq d(Hu, H(Hx_{2n+1})) + d(H(fx_{2n}), f(Hx_{2n})) + d(f(Hx_{2n}), fu).$

Thanks to (2.2) and (2.16)

(2.19) $$Hx_{2n} \to u, \quad fx_{2n} \to u$$

The pair $\{f, H\}$ is compatible, then

(2.20) $$d(H(fx_{2n}), f(Hx_{2n})) \to 0$$

Using the continuity of $f$ and (2.16), we have

(2.21) $$d(f(Hx_{2n}), fu) \to 0$$

Combining (2.17)-(2.20) together with (2.21) and letting $n \to +\infty$ in (2.18), we obtain

$$d(Hu, fu) \leq 0,$$

which means that $Hu = fu$. Again by the triangular inequality
(2.22)
$d(Hu, gu) \leq d(Hu, H(Hx_{2n+2})) + d(H(gx_{2n+1}), g(Gx_{2n+1})) + d(g(Gx_{2n+1}), gu).$

On the other hand

(2.23) $$Hx_{2n+1} \to u, \quad gx_{2n+1} \to u$$

The pair $\{g, H\}$ is compatible, then

(2.24) $$d(H(gx_{2n+1}), g(Gx_{2n+1})) \to 0$$

The continuity of $g$ together with (2.16) gives us

(2.25) $$d(g(Hx_{2n+1}), gu) \to 0$$

Combining (2.17)-(2.24) together with (2.25) and letting $n \to +\infty$ in (2.22), we obtain

$$d(Hu, gu) \leq 0,$$

which means that $Hu = gu$. We finish then by finding $fu = gu = Hu$, that is, $u$ is a coincidence point of $f$, $g$ and $H$. The proof of theorem 2.5 is proved. □

Now, we omit in the proof of theorem 2.5, the continuity of $f$, $g$ and $H$, and the compatibility of the pairs $\{f, H\}$ and $\{g, H\}$, and we replace them by other conditions in order to find the same result. This will be the purpose of the next theorem

**Theorem 2.6.** *Let $(X, \preceq)$ be a partially ordered set. Suppose that there exists a complete metric $d$ on $X$ such that $X$ is regular (i.e, if $\{z_n\}$ is a non-decreasing sequence in $X$ with respect to $\preceq$ such that $z_n \to z$ as $n \to +\infty$, then $z_n \preceq z$ for all $n \in \mathbb{N}$).*
*Let $f, g, H : X \to X$ be given mappings satisfying*
*(a) $fX \subseteq HX$, $gX \subseteq HX$,*
*(b) $HX$ is a closed subspace of $(X, d)$,*
*(c) $f$ and $g$ are weakly increasing with respect to $H$.*
*Suppose that for every $(x, y) \in X \times X$ such that $Hx$ and $Hy$ are comparable, we have (2.1) holds.*

*Then, $f$, $g$ and $H$ have a coincidence point $u \in X$, that is, $fu = gu = Hu$.*

*Proof.* We take the same sequences $\{x_n\}$ and $\{y_n\}$ as in the proof of theorem 2.5. In particular $\{Hx_n\}$ is a Cauchy sequence in the closed subspace $HX$, then there exists $v = Hu$,, $u \in X$ such that

(2.26) $$\lim_{n \to +\infty} Hx_n = v = Hu.$$

Thanks to (2.4), $\{Hx_n\}$ is non-decreasing sequence, then since it converges to $v = Hu$, we get

$$Hx_n \preceq Hu, \quad \forall\, n \in \mathbb{N},$$

so the terms $Hx_n$ and $Hu$ are comparable. Putting now $x = x_{2n}$ and $y = u$ in (2.1) and using (2.2) one can write

$$\begin{aligned} d(Hx_{2n+1}, gu) &= d(Tx_{2n}, gu) \\ &\leq \beta(d(Hx_{2n}, Hu))d(Hx_{2n}, Hu) \leq d(Hx_{2n}, Hu). \end{aligned}$$

Letting $n \to +\infty$ in the above inequality, using (2.26), we obtain

$$d(Hu, gu) = 0.$$

This means that $fu = Hu$. We use now the same strategy, one can write for $x = u$ and $y = x_{2n+1}$ in (2.1)

$$\begin{aligned} d(fu, Hx_{2n+2}) &= d(fu, gx_{2n+1}) \\ &\leq \beta(d(Hu, Hx_{2n+1}))d(Hu, Hx_{2n+1}) \\ &\leq d(Hu, Hx_{2n+1}). \end{aligned}$$

Similarly, we let $n \longrightarrow +\infty$ in the above inequality and we obtain

$$fu = Hu.$$

We conclude that $u$ is a coincidence point of $H$, $f$ and $g$, and then the proof of theorem 2.6 is completed. □

Now, we shall prove the existence and uniqueness theorem of a common fixed point.

**Theorem 2.7.** *In addition to the hypotheses of theorem 2.5, suppose that for any $(x, y) \in X \times X$, there exists $u \in X$ such that $fx \preceq fu$ and $fy \preceq fu$.*
*Then, $f$, $g$ and $H$ have a unique common fixed point, that is there exists a unique $z \in X$ such that $z = Hz = fz = gz$.*

*Proof.* Referring to theorem 2.5, the set of coincidence points is non-empty. We shall show that if $x^*$ and $y^*$ are coincidence points, that is, $Hx^* = fx^* = gx^*$ and $Hy^* = fy^* = gy^*$, then

(2.27) $$Hx^* = Hy^*.$$

By assumption, there exists $u_0 \in X$ such that

(2.28) $$fx^* \preceq fu_0, \quad fy^* \preceq fu_0.$$

Now, we proceed similar to the proof of theorem 2.5, we can immediately define the sequence $\{Hu_n\}$ as follows

(2.29) $$Hu_{2n+1} = fu_{2n}, \quad Hu_{2n+2} = gu_{2n+1}, \quad \forall\, n \in \mathbb{N}$$

Again, we have

(2.30) $$fx^* = Hx^* \preceq Hu_n, \quad fy^* = Hy^* \preceq Hu_n, \quad \forall\, n \in \mathbb{N}^*$$

Putting $x = u_{2n}$ and $y = x^*$ in (2.1) and using $\beta < 1$ and (2.30), we get
$$d(Hu_{2n+1}, Hx^*) = d(fu_{2n}, gx^*)) \leq \beta(d(Hu_{2n}, Hx^*))d(Hu_{2n}, Hx^*)$$
$$\leq d(Hu_{2n}, Hx^*).$$

This gives us

(2.31) $$d(Hu_{2n+1}, Hx^*) \leq d(Hu_{2n}, Hx^*)$$

Putting $x = x^*$ and $y = u_{2n}$ in (2.1), then similarly to the above, one can find

(2.32) $$d(Hu_{2n+2}, Hx^*) \leq d(Hu_{2n+1}, Hx^*).$$

Here we have used $Hx^* = gx^*$. We combine (2.31) to (2.32) to remark that

(2.33) $$d(Hu_{n+1}, Hx^*) \leq d(Hu_n, Hx^*).$$

Then the sequence $\{d(Hu_n, Hx^*)\}$ is non-increasing, so there exists $r \geq 0$ such that
$$d(Hu_n, Hx^*) \to r \quad \text{as} \quad n \to +\infty.$$

We know that $d(Hu_{2n+1}, Hx^*) \leq \beta(d(Hu_{2n}, Hx^*))d(Hu_{2n}, Hx^*)$, hence
$$\frac{d(Hu_{2n+1}, Hx^*)}{d(Hu_{2n}, Hx^*)} \leq \beta(d(Hu_{2n}, Hx^*)) < 1.$$

Letting $n \to +\infty$ in the above inequality, then we obtain
$$\lim_{n \to +\infty} \beta(d(Hu_{2n}, Hx^*)) = 1,$$

and since $\beta \in \mathcal{S}$, this implies that $d(Hu_{2n}, Hx^*) \to r = 0$. We then write

(2.34) $$d(Hu_n, Hx^*) \to 0 \quad \text{as} \quad n \to +\infty.$$

The same idea yields

(2.35) $$d(Hu_n, Hy^*) \to 0 \quad \text{as} \quad n \to +\infty.$$

(2.34)-(2.35) together with the fact that the limit is unique allows that (2.27) holds. Now, thanks to (2.29)-(2.34), we can write

(2.36) $$\lim_{n \to +\infty} fu_{2n} = Hx^* = Hy^*, \quad \lim_{n \to +\infty} gu_{2n+1} = Hx^* = Hy^*.$$

From the compatibility of the pairs $\{f, H\}$ and $\{g, H\}$, we obtain using (2.34)-(2.36)

(2.37) $$\lim_{n \to +\infty} d(H(fu_{2n}), f(Hu_{2n})) = 0, \quad \lim_{n \to +\infty} d(H(gu_{2n+1}), g(Hu_{2n+1})) = 0.$$

Let us denote
$$z = Hx^*.$$

By the triangular inequality, we have
$$d(Hz, fz) \leq d(Hz, H(fu_{2n})) + d(H(fu_{2n}), f(Hu_{2n})) + d(f(Hu_{2n}), fz),$$

Using (2.36)-(2.37) and the continuity of $f$ and letting $n \longrightarrow +\infty$ in the above inequality, we get
$$d(Hz, fz) \leq 0,$$

that is, $Hz = fz$ and $z$ is a coincidence point of $H$ and $f$. Again the triangular inequality gives us
$$d(Hz, gz) \leq d(Hz, H(gu_{2n+1})) + d(H(gu_{2n+1}), f(Hu_{2n+1})) + d(f(Hu_{2n+1}), gz),$$

Using (2.36)-(2.37) and the continuity of $g$ and letting $n \longrightarrow +\infty$ in the above inequality, we get
$$d(Hz, gz) \leq 0,$$

that is, $Hz = gz$ and $z$ is a coincidence point of $g$ and $H$. From (2.27), we have
$$z = Hx^* = Hz = fz = gz.$$

This proves that $z$ is a common fixed point of the mappings $H$, $g$ and $f$. Now our purpose is to check that such a point is unique. Suppose there is an another common fixed point $p$, that is

$$p = Hp = fp = gp.$$

This implies that $p$ is a coincidence point of $H$, $f$ and $g$. From (2.27), this implies that

$$Hp = Hz.$$

Hence, we get

$$p = Hp = Hz = z,$$

which yields the uniqueness of the common fixed point. The proof of theorem 2.7 is completed.
The next result is an immediate consequence of theorems 2.5 and 2.7 by taking $H = Id_X$. □

**Corollary 2.8.** *Let $(X, \preceq)$ be a partially ordered set. Suppose that there exists a complete metric $d$ on $X$. Let $f, g : X \to X$ be given mappings satisfying*
*(a) $f$, $g$ are continuous,*
*(b) $f$ and $g$ are weakly increasing.*
*Suppose that for every $(x, y) \in X \times X$ such that $x$ and $y$ are comparable, we have*

(2.38) $$d(fx, gy) \leq \beta(d(x,y)) \, d(x,y),$$

*where $\beta \in \mathcal{S}$. Suppose again that for each $x, y \in X$, there exists $z \in X$ which is comparable to $x$ and $y$. Then, $f$, $g$ have a unique common fixed point.*

**Remark 2.9.** *Taking $g = f$ in corollary 2.8, we find the result given in theorem 1.1.*

3. Application : existence for a common solution of integral equations

Consider the integral equations

(3.1) $$\begin{cases} u(t) = \int_0^T K_1(t, s, u(s))ds + h(t), & t \in [0, T] \\ u(t) = \int_0^T K_2(t, s, u(s))ds + h(t), & t \in [0, T] \end{cases}$$

where $T > 0$. The purpose of this section is to give an existence theorem for common solution of (3.1), using Corollary 2.8. This application is inspired in [2], [10]. Let consider the space $X = C(I)$ ($I = [0, T]$) of continuous functions defined on $I$. Obviously, this space with the metric given by

$$d(x, y) = \sup_{t \in I} |x(t) - y(t)|, \quad \forall \, x, y \in X,$$

is a complete metric space. $X = C(I)$ can also be equipped with the partial order $\preceq$ given by

$$\forall \, x, y \in X, \quad x \preceq y \iff x(t) \leq y(t), \quad \forall \, t \in I.$$

Moreover, in [12], it is proved that $(C(I), \preceq)$ is regular (the definition is given in theorem 2.6). Now we will state the following theorem which we find it in [11]

**Theorem 3.1.** *Suppose that the following hypotheses hold:*
*i) $K_1, K_2 : I \times I \times \mathbb{R} \to \mathbb{R}$ and $h : \mathbb{R} \to \mathbb{R}$ are continuous,*
*ii) for all $t, s \in I$, we have*

$$K_1(t, s, u(t)) \leq K_2\left(t, s, \int_0^T K_1(s, \tau, u(\tau))d\tau + h(s)\right).$$

$$K_2(t,s,u(t)) \leq K_1\left(t,s,\int_0^T K_2(s,\tau,u(\tau))d\tau + h(s)\right).$$

*iii)* there exists a continuous function $G : I \times I \to \mathbb{R}_+$ such that

$$|K_1(t,s,x) - K_2(t,s,y)| \leq G(t,s)\sqrt{\log[(x-y)^2 + 1]},$$

for all $t, s \in I$ and $x, y \in \mathbb{R}$ such that $y \preceq x$,
*iv)* $\sup_{t \in I} \int_0^T G^2(t,s)ds \leq \frac{1}{T}$.
Then the integral equations (3.1) have a solution $u^* \in C(I)$.

*Proof.* As mentioned above, the proof we find it in [11], with of course a leger modification at the end. We need to bring the quantities

$$fx(t) = \int_0^T K_1(t,s,u(s))ds + h(t), \quad t \in I.$$

$$gx(t) = \int_0^T K_2(t,s,u(s))ds + h(t), \quad t \in I.$$

Again, the authors [11] found that $f$ and $g$ are weakly increasing. Moreover, they obtained

$$d(fx, gy) \leq \sqrt{\log[d^2(x,y) + 1]}.$$

Now, the leger modification with respect to [11] is to choose the function $\beta$ as

$$\beta(t) = \frac{\sqrt{\log[t^2 + 1]}}{t}.$$

It is clear that with this choice, $\beta \in \mathcal{S}$. One can then write

$$d(fx, gy) \leq \beta(d(x,y))d(x,y).$$

Now, all the hypotheses of Corollary 2.8 are satisfied. Then, there exists $u^* \in C(I)$, a common fixed point of $f$ and $g$, that is $u^*$ is a solution to (3.1). □

H. Aydi

Université de Monastir, Institut supérieur d'informatique de Mahdia, Km 4, Réjiche, 5121 Mahdia, Tunisie

*E-mail address*: hassen.aydi@isima.rnu.tn